# THE THEORY OF CONNECTIONS AND THE PROBLEM OF EXISTENCE OF BACKLUND TRANSFORMATIONS FOR SECOND ORDER EVOLUTION EQUATIONS

## On the existence of Backlund Transformations for evolution equations

**A.K. Rybnikov**
arybnikov@mail.ru

Backlund transformations are used to search for solutions, particularly soliton solutions, of non-linear differential equations. In this paper we present an invariant geometrical theory of Backlund transformations for second order evolution equations. The main concept is that of connection defining the representation of zero curvature for a given partial differential equation. We systematically use the invariant analytical method of E. Cartan [C1,2] and G.F. Laptev [L1-5]. All differential-geometric considerations in this paper are local.

We show that the Backlund transformations of a second order evolution equation with one space variable

$$z_t - f(t,x,z,z_x,z_{xx}) = 0$$

can be represented in the following form

$$\left.\begin{aligned} y_t &= -y\cdot g^1_{01}(t,x,z,z_x) + \gamma^1_{00}(t,x,z,z_x)\,, \\ y_x &= -y\cdot\gamma^1_{11}(t,x,z,z_x) + \gamma^1_{01}(t,x,z,z_x) \end{aligned}\right\}\,,$$

where $g^1_{01}, \gamma^1_{00}, \gamma^1_{01}, \gamma^1_{11}$ are coefficients of a special connection defining the representation of zero curvature for the given evolution equation. Notice that Backlund transformations for evolution equations are more specific than Backlund transformations for general equations: $F(x^i, z, z_j, z_{kl}) = 0 \quad (i,j,\ldots = 1,2)$ (see Remark 1.1).

We prove that any second order evolution equation with one space variable, which admits Backlund transformations, is of the following form

$$z_t - L(t,x,z,z_x)\cdot z_{xx} - M(t,x,z,z_x) = 0\ .$$

In a special case, when the equation is of the form

$$z_t - z_{xx} - M(z,z_x) = 0\ ,$$

we state the problem of existence of standard Backlund transformations (for which $g^1_{01}, \gamma^1_{00}, \gamma^1_{01}, \gamma^1_{11}$ depend on $z$ and $z_x$ only). It is shown here that any equation of type $z_t - z_{xx} - M(z,z_x) = 0$, which admits standard Backlund transformations, has the following form

$$z_t - z_{xx} - \xi(z)\cdot(z_x)^2 - \eta(z)\cdot z_x - \zeta(z) = 0\ .$$

We prove that standard Backlund transformations for such equation exist if and only if the equation is of the form

$$z_t - z_{xx} - \frac{\eta''(z)}{\eta'(z)}\cdot(z_x)^2 - \eta(z)\cdot z_x = 0 \qquad (\eta'(z)\neq 0) \tag{Eq1}$$

or, of the form

$$z_t - z_{xx} - F'(z)\cdot(z_x)^2 - \eta(z)\cdot z_x - e^{-F(z)}\left(\int (A\eta(z)+B)e^{F(z)}dz + C\right) = 0, \tag{Eq2}$$

where $A,B,C = const$

It is possible to reduce (Eq1) to the Burgers equation

$$Z_t - Z_{xx} + Z\cdot Z_x = 0$$

by introducing a new dependent variable. Similarly, one can reduce (Eq2), in special case $\eta(z)\equiv k \quad (k = const)$, to the linear equation

$$Z_t - Z_{xx} - kZ_x - (kA+B)Z - C = 0\ .$$

These results have been partially announced in [R2, RS3].

## 1. Introduction.

Traditionally, the concept of Backlund map (and Backlund transformation) is interpreted in the following way. Suppose we are given a partial differential equation with $n$ independent variables $x^i \left(i, j, k, l, \ldots = 1, \ldots, n\right)$, their function $z$ and the partial derivatives $z_j, z_{kl}, \ldots$

$$F\left(x^i, z, z_j, z_{kl}, \ldots\right) = 0 . \qquad (1)$$

It is customary to define the Backlund map of the equation (1) as a system

$$\Psi_\xi \left(x^i, z, y, z_k, y_l, \ldots\right) = 0 , \qquad (2)$$

which is completely integrable over $y$ if and only if $z = z\left(x^1, \ldots, x^n\right)$ is a solution of the equation (1) (here $\xi \in U$, where $U$ is some set of indices). If any solution $y = y\left(x^1, \ldots, x^n\right)$ of the system (2) (where $z = z\left(x^1, \ldots, x^n\right)$ is any pre-assigned solution of the equation (1)) is at the same time a solution of a differential equation

$$\Phi\left(x^i, y, y_j, y_{kl}, \ldots\right) = 0 , \qquad (3)$$

then the Backlund map of the equation (1) is called the Backlund transformation of the equation (1) to the equation (3).

**Example.** Given the Burgers equation

$$z_t - z_{xx} + z \cdot z_x = 0 . \qquad \left(1'\right)$$

Consider the system

$$\left.\begin{aligned} y_t &= \frac{1}{2} y \left( z_x - \frac{1}{2} z^2 \right) + \frac{1}{2} c z , \\ y_x &= \frac{1}{2} y z - c \end{aligned}\right\} . \qquad \left(2'\right)$$

The system $\left(2'\right)$ is completely integrable over $y$ if and only if $z = z(t, x)$ is a solution of $\left(1'\right)$.

Using the second equation of the system $\left(2'\right)$, we get $z = \frac{2}{y}\left(y_x + c\right)$. Hence, $z_x = \frac{2}{y} y_{xx} - \frac{2}{y^2} y_x \left(y_x + c\right)$. Substituting the expressions for $z$ and $z_x$ into the first equation of the system $\left(2'\right)$, we get

$$y_t - y_{xx} + \frac{2}{y} \cdot \left(y_x\right)^2 + 2c \cdot \frac{y_x}{y} = 0 . \qquad \left(3'\right)$$

By choosing new dependent variable $Y = \frac{2c}{y}$, we obtain the Burgers equation

$$Y_t - Y_{xx} + Y \cdot Y_x = 0 . \qquad \square$$

The above-mentioned traditional definition of the Backlund map is not invariant with respect to non-degenerate coordinate changes; moreover, it seems to be artificial. It turns out that an invariant interpretation is possible in the differential-geometric setup. In this work we present the invariant theory of Backlund maps for second order evolution equations, i.e., for equations of the form

$$z_t - f\left(t, x^i, z, z_j, z_{kl}\right) = 0 \ . \qquad (4)$$

The main concept of this theory is the notion of connection defining the representation of zero curvature for the evolution equation (4). The theory of Backlund transformations is developed here as a special chapter in the theory of connections. We systematically employ the invariant analytical method of E. Cartan [C1,2] and G.F. Laptev [L1-5] and, in particular, Laptev's theory of structural forms of fibred spaces and connection forms (see [ELOS], or original works of Laptev [L1-5]). Nowadays this theory is the most powerful method of research in local differential geometry.

In 1979 F.A.E. Pirani, D.C. Robinson, and W.F. Shadwick [PRS] attempted to use the theory of connections for their investigation of Backlund transformations; however, their method was not so effective. In 1999, using the Cartan--Laptev method, we presented the invariant theory of Backlund maps for general (non-evolution) second order PDE [RS1,2]. The theory for the evolution equations presented in this work and the theory for general differential equations that we elaborated in 1999 have essential differences, although they are analogous in many respects.

While considering the evolution equation (4), we assume that $t, x^1, \ldots, x^n, z$ are adapted local coordinates on an $(n+2)$–dimensional general type bundle $E$ on a fibred $(n+1)$–dimensional base $M$ (the variables $t, x^1, \ldots x^n$ are local coordinates on the manifold $M$ ). We also assume that non-degenerate coordinate changes of the form

$$\tilde{t} = \varphi^0(t) \ ; \qquad \tilde{x}^i = \varphi^i\left(t, x^1, \ldots x^n\right) ; \qquad \tilde{z} = \varphi^{n+1}\left(t, x^1, \ldots, x^n, z\right)$$

are admissible transformations of local coordinates on $E$ . Notice that these coordinate transformations do not change the form (4) of the evolution equation.

We say that a section $\sigma \subset E$ is a solution of the evolution equation (4) if the equation (4) is identically satisfied on $\sigma$ . Notice that any section $\sigma \subset E$ can be defined by an equation of the form $z = z\left(t, x^1, \ldots x^n\right)$.

Let us briefly describe the structure of the paper. First, in Section 2, we recall the notion of *connection defining the representation of zero curvature* for a partial differential equation (PDE). Then, in Section 3, we consider the frame bundle of the first order $R^1 E$ on $E$ . We pick out two factor-manifolds of the manifold $R^1 E$: $R^* E$ and $\hat{R}^* E$ . Both $R^* E$ and $\hat{R}^* E$ are the principal bundles over a common base $J^* E$, where $J^* E$ is some factor-manifold of the 1-jet manifold $J^1 E$ . The Lie groups $\overline{GL}(n+1)$ and $\overline{SL}(n+1)$ are structural groups of $R^* R$ and $\hat{R}^* E$ respectively. The group $\overline{GL}(n+1)$ is a subgroup of the group $GL(n+1)$ . The group $\overline{SL}(n+1)$ is a subgroup of $\overline{GL}(n+1)$ and, at the same time, it is a subgroup of the group $SL(n+1)$ .

In Section 4 we study *special connections* in $R^* E$ and $\hat{R}^* E$ . These special connections generate corresponding connections in associated bundles $F\left(R^* E\right)$ and $F\left(\hat{R}^* E\right)$ with a typical fibre $F$ (and, in particular, in associated bundles with a one-dimensional typical fibre). Notice (see Lemma 3.1) that if $n = 1$, then the bundle $F\left(R^{*}E\right)$ $\left(\dim F = 1\right)$ does not exist; in this case only the bundle $F\left(\overset{\wedge}{R^{*}}E\right)$ $\left(\dim F = 1\right)$ exists.

In Section 5 we introduce the concept of *Backlund connection* corresponding to an evolution equation with one space variable

$$z_t - f\left(t,x,z,z_x,z_{xx}\right) = 0 \ .$$

A connection in associated bundle $F\left(\overset{\wedge}{R^{*}}E\right)$ $\left(\dim F = 1\right)$ is called a Backlund connection if it is generated by a special connection defining the representation of zero curvature.

The Pfaff equation

$$\underset{\sigma}{\tilde{\theta}} = 0 \ ,$$

where $\underset{\sigma}{\tilde{\theta}}$ is the connection form for a Backlund connection considered on a section $\sigma \subset E$ , is completely integrable if and only if the section $\sigma \subset E$ is a solution of the evolution equation (it is assumed that $\sigma \subset E$ is fixed). This Pfaff equation is equivalent (under a special choice of the principal forms) to a system of PDEs, which is traditionally called the "Backlund transformation" of (1). Yet, we prefer to use a slightly different terminology: we call this new system the *Backlund system.* We prove that the Backlund system for the second order evolution equation is of very special form (Theorem 1):

$$\left.\begin{aligned} y_t &= -y\cdot g^{1}_{01}\left(t,x,z,z_x\right) + \gamma^{1}_{00}\left(t,x,z,z_x\right) \ , \\ y_x &= -y\cdot\gamma^{1}_{11}\left(t,x,z,z_x\right) + \gamma^{1}_{01}\left(t,x,z,z_x\right) \end{aligned}\right\} \ ,$$

where $g^{1}_{01}, \gamma^{1}_{00}, \gamma^{1}_{01}, \gamma^{1}_{11}$ are coefficients of the special connection in $\overset{\wedge}{R^{*}}E$ $\left(n = 1\right)$ that define the representation of zero curvature (considered on the lifted section $\sigma^{1} \subset J^{1}E$ corresponding to the fixed section $\sigma \subset E$ ).

**Remark 1.1.** *It is known* [RS1,2] *that the Backlund system for a general PDE*

$$F\left(x^{i}, z, z_j, z_{kl}\right) = 0 \quad \left(i, j, \ldots = 1,2\right)$$

*can be written in the form*

$$y_i = y^{2}\,\Gamma^{1}_{2i} + y\Gamma_i - \Gamma^{2}_{1i} \quad \left(i = 1,2\right) \ ,$$

*where* $\Gamma_i, \Gamma^{1}_{2i}, \Gamma^{2}_{1i}$ *are the coefficients of a special connection defining the zero curvature representation for this PDE. By choosing new dependent variables in a convenient way one can transform this system to the system*

$$Y_i = \Gamma_i + C\,e^{Y}\Gamma^{1}_{2i} - \frac{1}{C}\cdot e^{-Y}\Gamma^{2}_{1i} \quad \left(i = 1,2\right) \ ,$$

*where* $C \neq 0$ *is a constant, or to the system*

$$X_i = \Gamma^{1}_{2i} - \Gamma^{2}_{1i} + \sin X \cdot \Gamma_i - \cos X \cdot \left(\Gamma^{1}_{2i} + \Gamma^{2}_{1i}\right) \quad \left(i = 1,2\right) \ .$$

We see that when a given PDE is an evolution PDE, the Backlund system is more special: in this case the differential structure on the manifold $E$ is different from the differential structure on $E$ in the general case.

In Section 6 we prove (Theorem 2) that any second order evolution equation with one space variable, that admits a Backlund transformation, is of the following form

$$z_t - L(t,x,z,z_x)\cdot z_{xx} - M(t,x,z,z_x) = 0 .$$

In Section 7 we consider a special case where the evolution equation is of the form

$$z_t - z_{xx} - M(z,z_x) = 0 .$$

For such equations we state the problem of existence of *standard Backlund transformations* (for which $g^1_{01}, \gamma^1_{00}, \gamma^1_{01}, \gamma^1_{11}$ depend on $z$ and $z_x$ only). It is proved (Theorem 3) that any equation of type $z_t - z_{xx} - M(z,z_x) = 0$, that admits standard Backlund transformations, has the following form

$$z_t - z_{xx} - \xi(z)\cdot(z_x)^2 - \eta(z)\cdot z_x - \zeta(z) = 0 .$$

Then, in Theorem 4, we prove that such equation admits standard Backlund transformations if and only if it is of the form (Eq1) or (Eq2). We notice (Remark 7.4) that it is possible to reduce (Eq1) to the Burgers equation $Z_t - Z_{xx} + Z\cdot Z_x = 0$. We also show (Remark 7.5) that (Eq2), in the special case $\eta(z) = k = cons$, can be reduced to the linear equation $Z_t - Z_{xx} - kZ_x - (kA + B)Z - C = 0$. The Backlund systems for (Eq1) and (Eq2) have the form (35) or the form (36) respectively (see Remark 7.6).

## 2. Connections in the associated bundles generated by connections defining representations of zero curvature.

Let $J^rE$ be the manifold of holonomic EMBEDr-jets of local sections of a manifold $E$, and let $B$ be a factor-manifold of $J^rE$ for some $r$ (or maybe $B = J^rE$). Recall that a connection in the principal bundle $P(B,G)$ is called *a connection defining the representation of zero curvature* for a given evolution equation (4) if its curvature forms vanish on the solutions $\sigma \subset E$ (more exactly, on the corresponding lifted sections of the solutions $\sigma \subset E$) of the equation (4) and only on them.

Any connection in the principal bundle $P(B,G)$ generates a connection in the associated bundle $F(P(B,G))$ with typical fibre $F$, where $F$ is the representation space of the structural group $G$. In particular, we consider generated connections in associated bundles with one-dimensional typical fibres.

Recall that a differential equation defining a 1-dimensional representation of Lie group $G$ has the form

$$dX - \xi_A(X)\cdot\overline{\omega}^A = 0 ,$$

where $\overline{\omega}^A$ are the invariant structural forms of the group $G$ $(A,B,\ldots = 1,\ldots,\dim G)$. Coefficients $\xi_A(X)$ satisfy Lie identities

$$\frac{d\xi_B}{dX}\cdot\xi_C - \frac{d\xi_C}{dX}\cdot\xi_B = \xi_A C^A_{BC} ; \tag{5}$$

here $C^A_{BC}$ are the structural constants of the Lie group $G$. The set of the fibred forms of $F(P(B,G))$ $(\dim F = 1)$ consists, as well-known (see for example [ELOS]), of a form

$$\theta = dY - \xi_A(Y)\cdot\omega^A .$$

Here $\omega^A$ are the fibred forms in $P(B,G)$, and the coefficients $\xi_A$ satisfy Lie identities (5). Let a connection in the associated bundle $F(P(B,G))$ $(\dim F = 1)$ is generated by a connection in the principal bundle $P(B,G)$. Then the set of the connection forms in $F(P(B,G))$ $(\dim F = 1)$ consists of a form

$$\tilde{\theta} = dY - \xi_A(Y)\cdot\tilde{\omega}^A ,$$

where $\tilde{\omega}^A$ are the connection forms in $P(B,G)$. The form $\tilde{\theta}$ satisfies the structural equation

$$d\tilde{\theta} = \tilde{\theta}\wedge\left(-\frac{d\xi_A}{dY}\cdot\tilde{\omega}^A\right) - \xi_A\Omega^A ,$$

where $\Omega^A$ are the curvature forms for the connection in $P(B,G)$.

**Remark 2.1.** *Suppose a connection in* $F(P(B,G))$ $(\dim F = 1)$ *is generated by a connection in* $P(B,G)$ *defining the representation of zero curvature for the evolution equation (4). Then Pfaff equation*

$$\tilde{\theta}_\sigma = 0$$

*considered on any section* $\sigma \subset E$ *is completely integrable if and only if the section* $\sigma \subset E$ *is a solution of equation (4) (it is understood that* $\sigma \subset E$ *is fixed).*

## 3. The principal bundles $R^*E$, $\hat{R}^*E$ and bundles associated with them.

Let us consider the sequence of jet bundles $J^1E$, $J^2E$, ... . Denote the adapted local coordinates of EMBED r-jet manifold $J^rE$ by $x^{\hat{i}}, z, p_{\hat{i}_1\ldots\hat{i}_\alpha}$, where $\alpha = 1,\ldots,r$ ($p_{\hat{i}_1\ldots\hat{i}_\alpha}$ are symmetric by subscripts). Here $\hat{i}, \hat{j},\ldots = 0,1,\ldots,n$ and $x^0 = t$. In particular, $x^{\hat{i}}, z, p_0, p_j$ are the local coordinates on 1-jet manifold $J^1E$. The same variables with the exception of $p_0$, i.e. variables $x^{\hat{i}}, z, p_j$, are the local coordinates of another manifold $\overset{*}{J^1}E$, which is a factor-manifold of $J^1E$.

For any section $\sigma \subset E$, defined by the equation $z = z(t, x^1,\ldots,x^n)$, one can consider the lifted sections $\sigma^r \subset J^rE$ $(r = 1,2,\ldots)$ defined by the equations

$$z = z(t, x^1,\ldots,x^n); \qquad p_{\hat{i}_1\ldots\hat{i}_\alpha} = z_{\hat{i}_1\ldots\hat{i}_\alpha} \quad (\alpha = 1,\ldots,r).$$

The following is a somewhat more general form of the equation (4) :

$$p_0 - f(t, x^i, z, p_j, p_{kl}) = 0 . \tag{6}$$

On the lifted section $\sigma^2 \subset J^2E$ of $\sigma \subset E$ the equation (6) has form (4).

Let $\omega^0, \omega^1,\ldots,\omega^n, \omega^{n+1}$ be the principal forms of the manifold $E$. They satisfy the structural equations

$$\left.\begin{array}{l} d\omega^0 = \omega^0 \wedge \omega^0_0 \,, \\ d\omega^i = \omega^0 \wedge \omega^i_0 + \omega^j \wedge \omega^i_j \,, \\ d\omega^{n+1} = \omega^0 \wedge \omega^{n+1}_0 + \omega^j \wedge \omega^{n+1}_j + \omega^{n+1} \wedge \omega^{n+1}_{n+1} \end{array}\right\} .$$

When we consider the regular prolongation of these equations (see [L3] to get acquainted with this procedure), a sequence of the structural forms of frame bundles of various orders appears. The set of structural forms of the first-order frame manifold $R^1E$ contains the forms $\omega^{\hat{i}}, \omega^{n+1}, \omega^{n+1}_{\hat{j}}$. These forms are also the principal forms of the 1-jet manifold $J^1E$. The same forms with the exception of $\omega^{n+1}_0$, i.e. the forms $\omega^{\hat{i}}, \omega^{n+1}, \omega^{n+1}_j$, are the principal forms of $\overset{*}{J^1} E$.

One can pick out factor-manifolds $R^*E$ and $\hat{R^*} E$ of the first-order frame manifold $R^1E$ so that:

- The forms $\omega^{\hat{i}}, \omega^{n+1}, \omega^{n+1}_j, \omega^0_0, \omega^i_0, \omega^i_j$ are the structural forms for the manifold $R^*E$;
- The forms $\omega^{\hat{i}}, \omega^{n+1}, \omega^{n+1}_j, \omega^i_0, \omega^i_j - \delta^i_j \omega^0_0$ are the structural forms for the manifold $\hat{R^*} E$.

Both $R^*E$ and $\hat{R^*} E$ are the principal bundles over the common base $\overset{*}{J^1} E$. If a point of the base $\overset{*}{J^1} E$ is fixed, the fibred forms of $R^*E$, i.e., the forms $\omega^0_0, \omega^i_0, \omega^i_j$, become the invariant structural forms of the Lie group $\overline{GL}(n+1)$ $\left(\overline{GL}(n+1) \subset GL(n+1)\right)$. The group $\overline{GL}(n+1)$ is the structural group of $R^*E$. Similarly, if a point of the base $\overset{*}{J^1} E$ is fixed, then the fibred forms of $\hat{R^*} E$, i.e., the forms $\omega^i_0, \omega^i_j - \delta^i_j \omega^0_0$, become the invariant structural forms of Lie group $\overline{SL}(n+1)$ ( $\overline{SL}(n+1) \subset \overline{GL}(n+1)$ and, at the same time, $\overline{SL}(n+1) \subset SL(n+1)$ ).

Let us consider the associated bundles $F(R^*E)$ $(\dim F = 1)$ and $F\left(\hat{R^*} E\right)$ $(\dim F = 1)$.

**LEMMA 3.1.** *One-dimensional representations of the Lie group* $\overline{GL}(1+1)$ *do not exist. So, if* $n=1$*, then the bundle* $F(R^*E)$ $(\dim F = 1)$ *does not exist.*

***PROOF:*** The proof is via *reductio ad absurdum*. Suppose that there is a 1-dimensional representation of $\overline{GL}(1+1)$. Recall that the structural equations of $\overline{GL}(1+1)$ are of the form

$$\left.\begin{aligned} &d\left(\overline{\omega}_1^1+\overline{\omega}_0^0\right)=0\ ,\\ &d\left(\overline{\omega}_1^1-\overline{\omega}_0^0\right)=0\ ,\\ &d\overline{\omega}_0^1=\overline{\omega}_0^1\wedge\left(\overline{\omega}_1^1-\overline{\omega}_0^0\right)\end{aligned}\right\}\ ;$$

note that the invariant structural forms of the group $\overline{GL}(1+1)$ (i.e., forms $\overline{\omega}_0^0\ ,\overline{\omega}_0^1\ ,\overline{\omega}_1^1$) can be replaced with the forms $\overline{\omega}_1^1+\overline{\omega}_0^0\ ,\ \overline{\omega}_1^1-\overline{\omega}_0^0\ ,\ \overline{\omega}_0^1$. The forms $\overline{\omega}_0^1\ ,\ \overline{\omega}_1^1-\overline{\omega}_0^0$ are at the same time the invariant structural forms of the group $SL(1+1)\subset\overline{GL}(1+1)$. Each of the forms $\overline{\omega}_1^1+\overline{\omega}_0^0\ ,\ \overline{\omega}_1^1-\overline{\omega}_0^0$ (or their linear combination with constant coefficients) can be considered as the invariant structural form of the 1-dimensional Lie group $G_1\subset\overline{GL}(1+1)$. The differential equation defining the 1-dimensional representation of $\overline{GL}(1+1)$ mus t be of the form

$$dY-\xi(Y)\cdot\overline{\omega}_0^1-\eta(Y)\cdot\left(\overline{\omega}_1^1-\overline{\omega}_0^0\right)-\zeta(Y)\cdot\left(\overline{\omega}_1^1+\overline{\omega}_0^0\right)=0\ . \qquad (7)$$

We claim that $\zeta(Y)\neq 0$ here. Really, if $\zeta(Y)=0$, then the differential equation (7) is of the following form

$$dY-\xi(Y)\cdot\overline{\omega}_0^1-\eta(Y)\cdot\left(\overline{\omega}_1^1-\overline{\omega}_0^0\right)=0\ ,$$

and this differential equation does not define the representation of 3-dimensional group $\overline{GL}(1+1)$ (it defines the representation of the 2-dimensional group $\overline{SL}(1+1)\subset\overline{GL}(1+1)$, if $\xi(Y)\neq 0$, or the representation of 1-dimensional group $G_1$, if $\xi(Y)=0$). EMBEDDATEdd/MM/yyyy HH:mm:ss

The coefficients $\xi(Y)$, $\eta(Y)$, $\zeta(Y)\neq 0$ must satisfy the Lie identities

$$\left.\begin{aligned} &\frac{d\xi}{dY}\cdot\zeta-\xi\cdot\frac{d\zeta}{dY}=0\ ,\\ &\frac{d\eta}{dY}\cdot\zeta-\eta\cdot\frac{d\zeta}{dY}=0\ ,\\ &\frac{d\eta}{dY}\cdot\xi-\eta\cdot\frac{d\xi}{dY}+\xi=0\end{aligned}\right\}.$$

One can rewrite these equations as follows

$$\left.\begin{aligned} &d\frac{\xi}{\zeta}=0\ ,\\ &d\frac{\eta}{\zeta}=0\ ,\\ &d\eta\cdot\xi-\eta\cdot d\xi=-\xi\cdot dY\end{aligned}\right\}. \qquad (8)$$

As $d\frac{\xi}{\zeta}=d\frac{\eta}{\zeta}=0$, we have

$$\xi=A\zeta\ ,\qquad \eta=B\zeta\qquad (A,B=const). \qquad (9)$$

Substituting (9) into the third equation of the system (8) we get, in light of $\zeta\neq 0$,

$$A=0\ .$$

Therefore, the differential equation (7) has the following form

$$dY-\zeta(Y)\cdot\left\{B\left(\overline{\omega}_1^1-\overline{\omega}_0^0\right)+A\left(\overline{\omega}_1^1+\overline{\omega}_0^0\right)\right\}=0\ ,$$

and it does not define a 1-dimensional representation of group $\overline{GL}(1+1)$: it defines a 1-dimension representation of group $G_1$.

This contradiction concludes the proof. □

**Remark 3.1.** *If* $n=1$, *one can represent fibred form in* $F\left(\overset{\wedge}{R^*}E\right)$ $(\dim F=1)$ *as follows*

$$\theta=\xi(Y)\cdot\left[dy+y\left(\omega_1^1-\omega_0^0\right)-\omega_0^1\right]\qquad(\xi(Y)\neq 0)\ , \tag{10}$$

*where* $y$ *is an antiderivative of the function* $-\frac{1}{\xi(Y)}$.

It is obvious that if $n=1$, the fibred form $\theta$ in $F\left(\overset{\wedge}{R^*}E\right)$ $(\dim F=1)$ has the following representation

$$\theta=dY-\xi(Y)\cdot\omega_0^1-\eta(Y)\cdot\left(\omega_1^1-\omega_0^0\right)\qquad(\xi(Y)\neq 0)\ .$$

In this case the Lie identities are as follows

$$\frac{dY}{\xi(Y)}=-d\frac{\eta(Y)}{\xi(Y)}\ .$$

Hence, (10) easily follows. □

## 4. Special connections in $R^*E$ and $\overset{\wedge}{R^*}E$. Connections in $F\left(\overset{\wedge}{R^*}E\right)$ $(\dim F=1)$ generated by special connections in $\overset{\wedge}{R^*}E$ $(n=1)$.

In $R^*E$ and $\overset{\wedge}{R^*}E$ one can consider *special connections*, i.e., connections for which all connection coefficients, except for coefficients of $\omega^0,\omega^1,\ldots,\omega^n$,**EMBED** are equal to zero. The class of special connections stands out invariantly [R1]. In the case $n=1$ the connection forms corresponding to the special connections in $R^*E$ (denoted by $\tilde{\omega}_0^0,\tilde{\omega}_0^1,\tilde{\omega}_1^1$) can be described by the system

$$\left.\begin{aligned}\tilde{\omega}_0^0&=\omega_0^0+\Gamma_{00}^0\omega^0\ ,\\ \tilde{\omega}_0^1&=\omega_0^1+\Gamma_{00}^1\omega^0+\Gamma_{01}^1\omega^1\ ,\\ \tilde{\omega}_1^1&=\omega_1^1+\Gamma_{01}^1\omega^0+\Gamma_{11}^1\omega^1\end{aligned}\right\}\ ,$$

and the connection forms corresponding to the special connections in $\overset{\wedge}{R^*}E$ (denoted by $\tilde{\tilde{\omega}}_0^1,\tilde{\tilde{\omega}}_1^1$) can be described by the system

$$\left.\begin{aligned}\tilde{\tilde{\omega}}_0^1&=\omega_0^1\qquad\quad+\gamma_{00}^1\omega^0+\gamma_{01}^1\omega^1\ ,\\ \tilde{\tilde{\omega}}_1^1&=\omega_1^1-\omega_0^0+g_{01}^1\omega^0+\gamma_{11}^1\omega^1\end{aligned}\right\}\ .$$

If $\omega^0=dt$, $\omega^1=dx$, one can assume, that $\omega_0^0=\omega_0^1=\omega_1^1=0$. In this case

$$\left.\begin{aligned}\tilde{\omega}_0^0 &= \Gamma_{00}^0\, dt,\\ \tilde{\omega}_0^1 &= \Gamma_{00}^1\, dt + \Gamma_{01}^1\, dx,\\ \tilde{\omega}_1^1 &= \Gamma_{01}^1\, dt + \Gamma_{11}^1\, dx\end{aligned}\right\};\qquad \left.\begin{aligned}\tilde{\tilde{\omega}}_0^1 &= \gamma_{00}^1\, dt + \gamma_{01}^1\, dx,\\ \tilde{\tilde{\omega}}_1^1 &= g_{01}^1\, dt + \gamma_{11}^1\, dx\end{aligned}\right\}. \qquad (11)$$

EMBED

**Remark 4.1.** *Let* $n=1$. *To define a special connection in* $R^*E$ *it is necessary and sufficient to define a special connection in* $\overset{\wedge}{R^*}E$ .

Really, it is easy to verify that if a special connection in $R^*E$ $(n=1)$ with the connection coefficients $\Gamma_{00}^0, \Gamma_{00}^1, \Gamma_{01}^1, \Gamma_{11}^1$ is given, then one can define the special connection in $\overset{\wedge}{R^*}E$ $(n=1)$ with the connection coefficients:

$$\gamma_{00}^1 = \Gamma_{00}^1, \qquad \gamma_{01}^1 = \Gamma_{01}^1,$$
$$g_{01}^1 = \Gamma_{01}^1 - \Gamma_{00}^0, \qquad \gamma_{11}^1 = \Gamma_{11}^1 .$$

Conversely, if a special connection in $\overset{\wedge}{R^*}E$ $(n=1)$ with the connection coefficients $g_{01}^1, \gamma_{00}^1, \gamma_{01}^1, \gamma_{11}^1$ is given, then one can define the special connection in $R^*E$ $(n=1)$ with the connection coefficients:

$$\Gamma_{00}^0 = \gamma_{01}^1 - g_{01}^1, \qquad \Gamma_{00}^1 = \gamma_{00}^1,$$
$$\Gamma_{01}^1 = \gamma_{01}^1, \qquad \Gamma_{11}^1 = \gamma_{11}^1 . \qquad \square$$

Special connections in principal bundles $R^*E$ and $\overset{\wedge}{R^*}E$ generate connections in associated bundles $F(R^*E)$ and $F\left(\overset{\wedge}{R^*}E\right)$ with the typical fibre $F$ (and, in particular, in the associated bundles with a 1-dimensional typical fibre). Notice that if $n=1$, the bundle $F(R^*E)$ $(\dim F=1)$ does not exist (see Lemma 3.1).

In this work we shall consider a connection in the associated bundle $F\left(\overset{\wedge}{R^*}E\right)$ $(\dim F=1)$ generated by a special connection in $\overset{\wedge}{R^*}E$ $(n=1)$. As follows from Remark 3.1, the curvature form corresponding to such connection can be written as

$$\tilde{\tilde{\theta}} = \xi(Y)\cdot\left(dy + y\tilde{\tilde{\omega}}_1^1 - \tilde{\tilde{\omega}}_0^1\right). \qquad (12)$$

Here, $y$ is an antiderivative of $-\dfrac{1}{\xi(Y)}$ $(\xi(Y)\neq 0)$. If $\omega^0 = dt$, $\omega^1 = dx$, $\omega_0^0 = \omega_0^1 = \omega_1^1 = 0$, then by (11) the form $\tilde{\tilde{\theta}}$ can be described as

$$\tilde{\tilde{\theta}} = \xi(Y)\cdot\left[dy - \left(-yg_{01}^1 + \gamma_{00}^1\right)dt - \left(-y\gamma_{11}^1 + \gamma_{01}^1\right)dx\right]. \qquad (13)$$

**Backlund connections and Backlund maps.**

Suppose a connection in the associated bundle $F\left(\overset{\wedge}{R^*} E\right)$ ($\dim F = 1$) is generated by a special connection in $\overset{\wedge}{R^*} E$ defining the representation of zero curvature for the evolution equation (6); then, the former connection in $F\left(\overset{\wedge}{R^*} E\right)$ ($\dim F = 1$) is called a *Backlund connection* of the evolution equation (6).

Notice (see Remark 2.1) that if $\tilde{\tilde{\theta}}$ is a connection form for the Backlund connection, then the Pfaff equation

$$\underset{\sigma}{\tilde{\tilde{\theta}}} = 0, \tag{14}$$

considered on a section $\sigma \subset E$, is completely integrable if and only if the section $\sigma \subset E$ is a solution of equation (6) (it is understood that $\sigma \subset E$ is fixed). The Pfaff equation (14) defines a mapping that takes each solution $\sigma \subset E$ of the equation (6) to the section $\underset{\sigma}{\Sigma} \subset F\left(\overset{\wedge}{R^*} E\right)$ ($\dim F = 1$), which is a solution of the Pfaff equation (14) considered on $\sigma \subset E$ (if $\sigma \subset E$ is fixed). We shall say that this mapping is the *Backlund map* corresponding to the evolution equation (6). We call the equation (14) the *Pfaff equation defining the Backlund map.*

Consider a Backlund connection corresponding to the evolution equation with one space variable:

$$p_0 - f(t, x, z, p_1, p_{11}) = 0 \,. \tag{15}$$

Let the principal forms of the bundle $\overset{\wedge}{R^*} E$ ($n = 1$) be contact forms

$$\omega^0 = dt; \;\; \omega^1 = dx; \;\; \omega^{1+1} = dz - p_0 dt - p_1 dx; \;\; \omega_1^{1+1} = dp_1 - p_{01} dt - p_{11} dx.$$

Then, by (13), the Pfaff equation (14) has the following form

$$dy - \left(-y g_{01}^1 + \gamma_{00}^1\right)dt - \left(-y\gamma_{11}^1 + \gamma_{01}^1\right)dx = 0 \,. \tag{16}$$

The coefficients $g_{01}^1, \gamma_{00}^1, \gamma_{01}^1, \gamma_{11}^1$ depend on $t, x, z, z_x$. The Pfaff equation (16) is equivalent to the system of PDEs

$$\left.\begin{aligned} y_t &= -y \cdot g_{01}^1(t, x, z, z_x) + \gamma_{00}^1(t, x, z, z_x)\,, \\ y_x &= -y \cdot \gamma_{11}^1(t, x, z, z_x) + \gamma_{01}^1(t, x, z, z_x) \end{aligned}\right\} . \tag{17}$$

The system (17) is exactly the system (2), which is traditionally called the "Backlund transformation", but our terminology is more invariant and, therefore, more natural.

So, the following theorem is valid

**THEOREM 1.** *The system (2) defining the Backlund map ( the Backlund system) for a second order evolution equation with one space variable can always be represented in the form (17), where* $g^1_{01}, \gamma^1_{00}, \gamma^1_{01}, \gamma^1_{11}$ *are the coefficients of a special connection in* $\hat{R}^* E$ $(n=1)$ *defining the representation of zero curvature for the evolution equation.*

## 6. Necessary conditions for existence of a Backlund map for a second order evolution equation with one space variable.

Let us consider second order evolution equations with one space variable, i.e., equations of the form (15). Ask the following question. What evolution equations of form (15) admit connections defining the representations of zero curvature (and therefore admit Backlund maps)?

Let us examine the structural equations of a special connection in $\hat{R}^* E$ $(n=1)$ (we assume that principal forms are contact). These equations can be written in following form (dots denote the sum of terms that contain the external products of the principal forms different from $dt\wedge dx$)

$$\left.\begin{aligned} d\tilde{\omega}^1_0 - \tilde{\omega}^1_0\wedge\tilde{\omega}^1_1 &= 2\rho^1_{001}\, dt\wedge dx + \dots ,\\ d\tilde{\omega}^1_1 &= 2\rho_{01}\, dt\wedge dx + \dots \end{aligned}\right\}, \tag{18}$$

where

$$\left.\begin{aligned} 2\rho_{01} &= \frac{\partial\gamma^1_{11}}{\partial t} - \frac{\partial g^1_{01}}{\partial x} + \frac{\partial\gamma^1_{11}}{\partial z}\cdot p_0 - \frac{\partial g^1_{01}}{\partial z}\cdot p_1 + \frac{\partial\gamma^1_{11}}{\partial p_1}\cdot p_{01} - \frac{\partial g^1_{01}}{\partial p_1}\cdot p_{11} ,\\ 2\rho^1_{001} &= \frac{\partial\gamma^1_{01}}{\partial t} - \frac{\partial\gamma^1_{00}}{\partial x} + \frac{\partial\gamma^1_{01}}{\partial z}\cdot p_0 - \frac{\partial\gamma^1_{00}}{\partial z}\cdot p_1 + \frac{\partial\gamma^1_{01}}{\partial p_1}\cdot p_{01} - \frac{\partial\gamma^1_{00}}{\partial p_1}\cdot p_{11} -\\ &\quad -\gamma^1_{00}\cdot\gamma^1_{11} + \gamma^1_{01}\cdot g^1_{01} \end{aligned}\right\}. \tag{19}$$

Notice that $\rho_{01}$ and $\rho^1_{001}$ are the components of a tensor (this tensor is a subtensor of the curvature tensor).

It is well known that if the forms $\omega^{\hat{i}}, \omega^{n+1}, \omega^{n+1}_{\hat{j}}, \dots, \omega^{n+1}_{\hat{i}_1\dots\hat{i}_k}$ are contact, then the lifted sections $\sigma^{k-1}\subset J^{k-1}E$ of the sections $\sigma\subset E$ (and only the sections $\sigma^{k-1}\subset J^{k-1}E$) are the integral manifolds of the Pfaff system

$$\omega^{n+1} = \omega^{n+1}_{\hat{j}} = \dots = \omega^{n+1}_{\hat{i}_1\dots\hat{i}_k} = 0 .$$

So, on the lifted section of any section $\sigma\subset E$ the equations (18) have the form

$$\left.\begin{aligned} d\tilde{\omega}^1_0 - \tilde{\omega}^1_0\wedge\tilde{\omega}^1_1 &= 2\underset{\sigma}{\rho^1_{001}}\, dt\wedge dx,\\ d\tilde{\omega}^1_1 &= 2\underset{\sigma}{\rho_{01}}\, dt\wedge dx \end{aligned}\right\}.$$

Here the coefficients $2\underset{\sigma}{\rho^1_{001}}$, $2\underset{\sigma}{\rho_{01}}$ satisfy the equations (19) under the condition

$$p_0 = z_t , \qquad p_1 = z_x , \qquad p_{01} = z_{tx} , \qquad p_{11} = z_{xx} . \tag{20}$$

Curvature forms vanish on the section $\sigma \subset E$ if and only if

$$\underset{\sigma}{\rho_{01}} = 0\,, \qquad \underset{\sigma}{\rho^{1}_{001}} = 0\,. \tag{21}$$

The system (21) is a system of the form

$$\left.\begin{aligned} &\left.\frac{\partial \gamma^{1}_{11}}{\partial z}\right|_{p_1 = z_x} \cdot z_t + \left.\frac{\partial \gamma^{1}_{11}}{\partial p_1}\right|_{p = z_x} \cdot z_{tx} - \left.\frac{\partial g^{1}_{01}}{\partial p_1}\right|_{p_1 = z_x} \cdot z_{xx} + U(t,x,z,z_x) = 0\,, \\ &\left.\frac{\partial \gamma^{1}_{01}}{\partial z}\right|_{p_1 = z_x} \cdot z_t + \left.\frac{\partial \gamma^{1}_{01}}{\partial p_1}\right|_{p_1 = z_x} \cdot z_{tx} - \left.\frac{\partial \gamma^{1}_{00}}{\partial p_1}\right|_{p_1 = z_x} \cdot z_{xx} + V(t,x,z,z_x) = 0 \end{aligned}\right\}. \tag{22}$$

The special connection in $\overset{\wedge}{R^{*} E}$ $(n=1)$ defines a representation of zero curvature for the equation (15) if and only if the system (21) is equivalent to the equation (15). This equivalence holds if and only if we can get the first equation of the system (22) by multiplying the equation (15) by $\frac{\partial \gamma^{1}_{11}}{\partial z}$, and get the second equation of the system (22) by multiplying (15) by $\frac{\partial \gamma^{1}_{01}}{\partial z}$. So, the equation (15) must be of the form

$$z_t - L(t,x,z,z_x) \cdot z_{xx} - M(t,x,z,z_x) = 0, \tag{23}$$

and the coefficients $g^{1}_{01}, \gamma^{1}_{00}, \gamma^{1}_{01}, \gamma^{1}_{11}$ must form a solution of the PDE system

$$\left.\begin{aligned} &\frac{\partial \gamma^{1}_{11}}{\partial p_1} = \frac{\partial \gamma^{1}_{01}}{\partial p_1} = 0\,, \\ &L \cdot \frac{\partial \gamma^{1}_{11}}{\partial z} = \frac{\partial g^{1}_{01}}{\partial p_1}\,, \\ &M \cdot \frac{\partial \gamma^{1}_{11}}{\partial z} = \frac{\partial g^{1}_{01}}{\partial z} \cdot p_1 - \frac{\partial \gamma^{1}_{11}}{\partial t} + \frac{\partial g^{1}_{01}}{\partial x}\,, \\ &L \cdot \frac{\partial \gamma^{1}_{01}}{\partial z} = \frac{\partial \gamma^{1}_{00}}{\partial p_1}\,, \\ &M \cdot \frac{\partial \gamma^{1}_{01}}{\partial z} = \frac{\partial \gamma^{1}_{00}}{\partial z} \cdot p_1 - \frac{\partial \gamma^{1}_{01}}{\partial t} + \frac{\partial \gamma^{1}_{00}}{\partial x} + \gamma^{1}_{00} \cdot \gamma^{1}_{11} - \gamma^{1}_{01} \cdot g^{1}_{01} \end{aligned}\right\} \tag{24}$$

under the additional condition

$$\left(\frac{\partial \gamma^{1}_{11}}{\partial z}\right)^2 + \left(\frac{\partial \gamma^{1}_{01}}{\partial z}\right)^2 \neq 0\,. \tag{25}$$

Thus, the following theorem holds:

**THEOREM 2.** *Any second order evolution equation with one space variable, which admits a Backlund map, is of the form*

$$z_t - L(t,x,z,z_x) \cdot z_{xx} - M(t,x,z,z_x) = 0\,,$$

*and the Backlund maps exist if and only if the system (24) has a solution under the condition (25).*

**Remark 6.1.** *Notice the coefficients* $\gamma^{1}_{01}$ *and* $\gamma^{1}_{11}$ *do not depend on* $p_1$ (*see the first equations of the system (24)*).

## 7. The problem of existence of Backlund transformations for evolution equations of the form $z_t - L(z,z_x)\cdot z_{xx} - M(z,z_x) = 0$.

Suppose we are given a special connection in $\overset{\wedge}{R^{*}E}$ $(n=1)$ defining the representation of zero curvature for the evolution equation (23). If the coefficients $g^{1}_{01}, \gamma^{1}_{00}, \gamma^{1}_{01}, \gamma^{1}_{11}$ depend only on $z$ and $p_1$, we say that this connection is *standard*.

**Remark 7.1.** *If the special connection is standard, then the coefficients* $\gamma^{1}_{11}$ *and* $\gamma^{1}_{01}$ *are the functions of* $z$ *only, as they do not depend on* $p_1$ *(see Remark 6.1).*

**Remark 7.2.** *If there exists a standard connection for the equation (23), then this equation is of the form*

$$z_t - L(z,z_x)\cdot z_{xx} - M(z,z_x) = 0 . \tag{26}$$

It is warranted by the form of the system (24). □

We say that a *Backlund map is standard* if it corresponds to a standard connection. A standard Backlund map is defined by a Backlund system of the form

$$\left.\begin{aligned} y_t &= -y\cdot g^{1}_{01}(z,z_x) + \gamma^{1}_{00}(z,z_x) , \\ y_x &= -y\cdot \gamma^{1}_{11}(z) + \gamma^{1}_{01}(z) \end{aligned}\right\}, \tag{27}$$

under the additional condition

$$\left(\frac{d\gamma^{1}_{11}}{dz}\right)^2 + \left(\frac{d\gamma^{1}_{01}}{dz}\right)^2 \neq 0 . \tag{25'}$$

**Remark 7.3.** *A standard Backlund map is always a Backlund transformation.*

It follows from $(25')$ and the second equation of the system (27) that $z = z(y,y_x)$. Hence $z_x = z_x(y,y_x,y_{xx})$. Substituting the expressions for $z$ and $z_x$ into the first equation of the system (27), we obtain a second order evolution equation in $y$. □

In this work we solve the problem of existence of the standard Backlund transformations for the evolution equation (26) in the special case of $L(z,p_1)\equiv 1$, i.e., for the equations

$$z_t - z_{xx} - M(z,z_x) = 0 . \tag{28}$$

Let us first prove the following theorem:

**THEOREM 3.** *Any equation of type (28), that admits a standard Backlund transformation, is of the form*

$$z_t - z_{xx} - \xi(z)\cdot(z_x)^2 - \eta(z)\cdot z_x - \zeta(z) = 0 . \tag{29}$$

*Moreover, if* $\frac{d\gamma^1_{11}}{dz} \neq 0$ *, then* $\zeta(z) = 0$.

***PROOF:*** By Theorem 2 (and Remarks 7.1 and 7.3), an evolution equation of type (28) admits a standard Backlund transformation if and only if the PDE system (24), where $L = 1$ and $M = M(z, p_1)$, has a solution for which $g^1_{01} = g^1_{01}(z, p_1)$, $\gamma^1_{00} = \gamma^1_{00}(z, p_1)$, $\gamma^1_{01} = \gamma^1_{01}(z)$, $\gamma^1_{11} = \gamma^1_{11}(z)$, and, moreover, the condition ($25'$) holds. This PDE system has the following form

$$\left.\begin{aligned} &\frac{d\gamma^1_{11}}{dz} = \frac{\partial g^1_{01}}{\partial p_1}, \\ &M \cdot \frac{d\gamma^1_{11}}{dz} = \frac{\partial g^1_{01}}{\partial z} \cdot p_1, \\ &\frac{d\gamma^1_{01}}{dz} = \frac{\partial \gamma^1_{00}}{\partial p_1}, \\ &M \cdot \frac{d\gamma^1_{01}}{dz} = \frac{\partial \gamma^1_{00}}{\partial z} \cdot p_1 + \gamma^1_{00} \cdot \gamma^1_{11} - \gamma^1_{01} \cdot g^1_{01} \end{aligned}\right\}. \tag{30}$$

If the system (30) has a solution, then it is possible that for this solution $\frac{d\gamma^1_{11}}{dz} \neq 0$ (the case 1) or $\frac{d\gamma^1_{11}}{dz} = 0$ (the case 2). Let us treat these cases separately.

**Case 1.** Let be $\frac{d\gamma^1_{11}}{dz} \neq 0$; then we have (see the first equation of the system (30))

$$g^1_{01} = \frac{d\gamma^1_{11}}{dz} \cdot p_1 + \varphi(z). \tag{31}$$

Here $\varphi(z)$ is some function.

Substituting (31) into the second equation of the system (30), we get

$$M = \left(\frac{d\gamma^1_{11}}{dz}\right)^{-1} \cdot \frac{d^2\gamma^1_{11}}{dz^2} \cdot (p_1)^2 + \left(\frac{d\gamma^1_{11}}{dz}\right)^{-1} \cdot \frac{d\varphi}{dz} \cdot p_1 .$$

Hence, in the Case 1

$$M = \xi(z) \cdot (p_1)^2 + \eta(z) \cdot p_1 ,$$

where $\xi(z) = \left(\frac{d\gamma^1_{11}}{dz}\right)^{-1} \cdot \frac{d^2\gamma^1_{11}}{dz^2}$, $\eta(z) = \left(\frac{d\gamma^1_{11}}{dz}\right)^{-1} \cdot \frac{d\varphi}{dz}$.

**Case 2.** Let $\frac{d\gamma^1_{11}}{dz} = 0$; then it follows from ($25'$) that $\frac{d\gamma^1_{01}}{dz} \neq 0$.

In this case

$$\gamma^1_{11} = a = cons, \tag{32}$$

and, moreover (see the first and the second equations of the system (30)),

$$\frac{\partial g^1_{01}}{\partial z} = \frac{\partial g^1_{01}}{\partial p_1} = 0 .$$

Therefore,

$$g^{1}_{01} = b = const \tag{33}$$

From the third equation of the system (30) it follows that

$$\gamma^{1}_{00} = \frac{d\gamma^{1}_{01}}{dz} \cdot p_1 + \psi(z), \tag{34}$$

where $\psi(z)$ is some function.

Substituting (32), (33), and (34) into the fourth equation of the system (30), we get

$$M = \left(\frac{d\gamma^{1}_{01}}{dz}\right)^{-1} \cdot \frac{d^2\gamma^{1}_{01}}{dz^2} \cdot (p_1)^2 + \left(\left(\frac{d\gamma^{1}_{01}}{dz}\right)^{-1} \cdot \frac{d\psi}{dz} + a\right) \cdot p_1 + \left(\frac{d\gamma^{1}_{01}}{dz}\right)^{-1} \cdot \left(a\psi - b\gamma^{1}_{01}\right).$$

Hence, in the Case 2

$$M = \xi(z) \cdot (p_1)^2 + \eta(z) \cdot p_1 + \zeta(z),$$

where $\xi(z) = \left(\frac{d\gamma^{1}_{01}}{dz}\right)^{-1} \cdot \frac{d^2\gamma^{1}_{01}}{dz^2}$, $\eta(z) = \left(\frac{d\gamma^{1}_{01}}{dz}\right)^{-1} \cdot \frac{d\psi}{dz} + a$,

$\zeta(z) = \left(\frac{d\gamma^{1}_{01}}{dz}\right)^{-1} \cdot \left(a\psi - b\gamma^{1}_{01}\right)$.

We therefore see that the Theorem 3 is valid. □

We are interested in the following two equations

$$z_t - z_{xx} - \frac{\eta''(z)}{\eta'(z)} \cdot (z_x)^2 - \eta(z) \cdot z_x = 0 \qquad \left(\eta'(z) \neq 0\right), \tag{Eq1}$$

$$z_t - z_{xx} - F'(z) \cdot (z_x)^2 - \eta(z) \cdot z_x - e^{-F(z)}\left(\int \left(A\eta(z) + B\right) e^{F(z)} dz + C\right) = 0, \tag{Eq2}$$

where $A, B, C = cons$; here $\int \left(A\eta(z) + B\right) e^{F(z)} dz$ is an arbitrary chosen antiderivative.

**Remark 7.4.** *Equation (Eq1) can be reduced to the Burgers equation ($1'$) in $Z$,*

$$Z_t - Z_{xx} + Z \cdot Z_x = 0$$

*by the substitution $Z = -\eta(z)$.*

**Remark 7.5.** *In the special case where $\eta(z) \equiv k$ the equation (Eq2) has the following form*

$$z_t - z_{xx} - F'(z) \cdot (z_x)^2 - k z_x - e^{-F(z)}\left((kA + B)\int e^{F(z)} dz + C\right) = 0. \tag{Eq 2'}$$

*The equation (Eq $2'$) is reduced to the linear equation in Z,*

$$Z_t - Z_{xx} - k Z_x - (kA + B) Z - C = 0$$

*by the substitution $Z = \int e^{F(z)} dz$; here $\int e^{F(z)} dz$ is an arbitrary chosen antiderivative.*

Now we can state

**THEOREM 4.** *An evolution equation of type*

$$z_t - z_{xx} - M(z, z_x) = 0$$

*admits a standard Backlund transformation if and only if it is of the form (Eq1) or, of the form (Eq2).*

***PROOF:*** The proof follows from Theorem 3 and Lemmas A and B (see Appendices A and B). □

**Remark 7.6.** *For (Eq1) the Backlund system has the following form*

$$\left.\begin{aligned} y_t &= \frac{1}{2}(-y+c_1)\cdot\left[\eta'(z)\cdot z_x+\frac{1}{2}\left(\eta^2(z)-\left(\frac{c_2}{a}\right)^2\right)\right]+\frac{a}{2}\left(\eta(z)+\frac{c_2}{a}\right), \\ y_x &= \frac{1}{2}(-y+c_1)\cdot\left(\eta(z)-\frac{c_2}{a}\right)+a \end{aligned}\right\}, \tag{35}$$

*where* $a, c_1, c_2$ *are constants (* $a \neq 0$*).*

*For (Eq2) the Backlund system has the form*

$$\left.\begin{aligned} y_t &= \left[\left(A^2+B\right)y+c\left(\int\left(\eta(z)-A\right)e^{F(z)}dz+c_2\right)\right]+ce^{F(z)}\cdot z_x, \\ y_x &= -Ay+c\left(\int e^{F(z)}dz+c_1\right) \end{aligned}\right\}, \tag{36}$$

*where* $c, c_1, c_2$ *are constants such that* $c \neq 0$, *and* $\tilde{n}_1$, $\tilde{n}_2$ *satisfy the relation* $\left(A^2+B\right)c_1+Ac_2=C$.

By using (17) and (a10), (a11) from Appendix A, we get (35). Similarly, using (17) and (b11) from Appendix B we obtain (36). □

**Remark 7.7.** *The transformation (35) takes the equation (Eq1) to the equation*

$$y_t-y_{xx}-y_x\cdot\left(\frac{c_2}{a}-2\frac{y_x-a}{y-c_1}\right)=0,$$

*which is reducible to the Burgers equation (*$1'$*) in* $Z$ *by the substitution*

$$Z=-\left(\frac{2a}{y-c_1}+\frac{c_2}{a}\right).$$

**Remark 7.8.** *The transformation (36) takes the equation (Eq* $2'$*) to the linear equation*

$$y_t-y_{xx}-ky_x-\left(kA+B\right)y-cc_2=0.$$

In appendices A and B we consider the problem of existence of a standard connection corresponding to a given evolution equation of type (29). We investigate separately the problem of existence a standard connections, for which $\frac{d\gamma^1_{11}}{dz} \neq 0$ (in appendix A) and $\frac{d\gamma^1_{11}}{dz} = 0$ (in appendix B).

**APPENDIX A. Standard connections corresponding to the evolution equations of type (29) and satisfying $\frac{d\gamma^1_{11}}{dz} \neq 0$.**

EMBED Recall (see Theorem 3) that if $\frac{d\gamma^1_{11}}{dz} \neq 0$, then the evolution equation is of the following form

$$z_t - z_{xx} - \xi(z)\cdot(z_x)^2 - \eta(z)\cdot z_x = 0\,. \qquad \text{(a1)}$$

Recall also that the required standard connection exists if and only if the PDE system (30) has a solution (under the condition $M = \xi(z)\cdot(p_1)^2 + \eta(z)\cdot p_1$), for which $\frac{d\gamma^1_{11}}{dz} \neq 0$ (see Theorem 2 and Remarks 7.1 and 7.2). Let us write up this PDE system:

$$\left.\begin{aligned}
&(\gamma^1_{11})' = \frac{\partial g^1_{01}}{\partial p_1}\,,\\
&(\xi\cdot p_1 + \eta)\cdot(\gamma^1_{11})' = \frac{\partial g^1_{01}}{\partial z}\,,\\
&(\gamma^1_{01})' = \frac{\partial \gamma^1_{00}}{\partial p_1}\,,\\
&(\xi\cdot(p_1)^2 + \eta\cdot p_1)\cdot(\gamma^1_{01})' = \frac{\partial \gamma^1_{00}}{\partial z}\cdot p_1 + \gamma^1_{00}\cdot\gamma^1_{11} - \gamma^1_{01}\cdot g^1_{01}
\end{aligned}\right\}; \qquad \text{(a2)}$$

here $\frac{d\gamma^1_{11}}{dz}$, $\frac{d\gamma^1_{01}}{dz}$ are denoted by $(\gamma^1_{11})'$, $(\gamma^1_{01})'$ respectively.

It follows from the first and the third equations of the system (a2) that

$$\left.\begin{aligned}
g^1_{01} &= (\gamma^1_{11})'\cdot p_1 + \varphi(z),\\
\gamma^1_{00} &= (\gamma^1_{01})'\cdot p_1 + \psi(z)
\end{aligned}\right\},$$

where $\varphi(z)$ and $\psi(z)$ are certain functions.

Substituting these expressions into the 2nd the 4-th equations of the system (a2), we get

$$\left.\begin{aligned}
&\xi\cdot(\gamma_{11}^1)'\cdot p_1+\eta\cdot(\gamma_{11}^1)'=(\gamma_{11}^1)''\cdot p_1+\varphi',\\
&\xi\cdot(\gamma_{01}^1)'\cdot(p_1)^2+\eta\cdot(\gamma_{01}^1)'\cdot p_1=(\gamma_{01}^1)''\cdot(p_1)^2+\left(\psi'+\gamma_{11}^1\cdot(\gamma_{01}^1)'-\gamma_{01}^1\cdot(\gamma_{11}^1)'\right)\cdot p_1+\\
&+\gamma_{11}^1\cdot\psi-\gamma_{01}^1\cdot\varphi
\end{aligned}\right\}.$$

Both sides in these relations are polynomials in $p_1$. Compare the coefficients on both sides:

$$\left.\begin{aligned}
&\xi\cdot(\gamma_{01}^1)'=(\gamma_{01}^1)'',\\
&\xi\cdot(\gamma_{11}^1)'=(\gamma_{11}^1)'',\\
&\eta\cdot(\gamma_{01}^1)'=\psi'+\gamma_{11}^1\cdot(\gamma_{01}^1)'-\gamma_{01}^1\cdot(\gamma_{11}^1)',\\
&\eta\cdot(\gamma_{11}^1)'=\varphi',\\
&\gamma_{11}^1\cdot\psi-\gamma_{01}^1\cdot\varphi=0
\end{aligned}\right\}. \qquad \text{(a3)}$$

**Remark A1.** *Notice that* $(\gamma_{01}^1)'$ *and* $(\gamma_{11}^1)'$ *are proportional:*

$$(\gamma_{01}^1)'=c_1(\gamma_{11}^1)', \qquad \text{(a4)}$$

*where* $c_1=const$ *and, therefore,*

$$\left.\begin{aligned}
&g_{01}^1=(\gamma_{11}^1)'\cdot p_1+\varphi(z),\\
&\gamma_{00}^1=c_1(\gamma_{11}^1)'\cdot p_1+\psi(z)
\end{aligned}\right\}. \qquad \text{(a5)}$$

***PROOF:*** This statement follows from the 1st and 2nd equations of (a3). In case of $(\gamma_{01}^1)'=0$ our assertion is evident. Therefore, let $(\gamma_{01}^1)'\neq 0$, and then we have

$$\frac{(\gamma_{01}^1)''}{(\gamma_{01}^1)'}=\frac{(\gamma_{11}^1)''}{(\gamma_{11}^1)'}.$$

Hence,

$$\left(\ln\left|(\gamma_{01}^1)'\right|\right)'=\left(\ln\left|(\gamma_{11}^1)'\right|\right)',$$

and (a5) holds. □

**Remark A2.** *It follows from (a4) that*

$$\gamma_{01}^1=c_1\gamma_{11}^1+a, \qquad \text{(a6)}$$

*where* $a$, $c_1=cons$, *and, moreover,*

$$a\neq 0. \qquad \text{(a7)}$$

***PROOF:*** Since (a6) is evident, it is sufficient to prove (a7). The proof is by *reductio ad absurdum.* Assume the converse, i.e., let $\tilde{a}=0$. Then

$$\gamma_{01}^1=c_1\gamma_{11}^1, \qquad \text{(a8)}$$

and, in view of $\gamma_{11}^1\neq 0$, the last equation of (a3)EMBED implies that

$$\psi = c_1\varphi \ .$$

Therefore, by invoking (a5), we have

$$\gamma^{\,1}_{\,00} = c_1\, g^1_{01} \ . \tag{a9}$$

If $\omega^0 = dt,\ \ \omega^1 = dx$, then by using (a8), (a9), and (11), we get

$$\tilde{\tilde{\omega}}^1_0 = c_1\tilde{\tilde{\omega}}^1_1 \ ,$$

The last relation always holds ( not only when $\omega^0 = dt,\ \ \omega^1 = dx$).

So, if $a = 0$, then the connection forms $\tilde{\tilde{\omega}}^1_0$ , $\tilde{\tilde{\omega}}^1_1$ are linearly dependent. But, in general, they are independent. This contradiction concludes the proof. □

**Remark A3.** *In light of Remarks A1 and A2, it follows from (a3) that*

$$\gamma^{\,1}_{\,11} = \frac{1}{2}\left(\eta(z) - \frac{c_2}{a}\right) \qquad (c_2 = const), \tag{a10}$$

*and, therefore, since* $\left(\gamma^{\,1}_{\,11}\right)' \neq 0$, *we have*

$\eta'(z) \neq 0$ .

*Moreover,*

$$\left.\begin{aligned} \gamma^{\,1}_{\,01} &= \frac{c_1}{2}\left(\eta(z) - \frac{c_2}{a}\right) + a \ , \\ \gamma^{\,1}_{\,00} &= \frac{c_1}{2}\left[\eta'(z)\cdot p_1 + \frac{1}{2}\left(\eta^2(z) - \left(\frac{c_2}{a}\right)^2\right)\right] + \frac{a}{2}\left(\eta(z) + \frac{c_2}{a}\right) , \\ g^1_{01} &= \frac{1}{2}\cdot\eta'(z)\cdot p_1 + \frac{1}{4}\left(\eta^2(z) - \left(\frac{c_2}{a}\right)^2\right) \end{aligned}\right\} . \tag{a11}$$

***PROOF*:** Substituting (a6) into the 3rd equation of the system (a3), and taking into account the 4th equation of (a3), we get

$$\left(\psi - a\gamma^{\,1}_{\,11} - c_1\varphi\right)' = 0 \ .$$

Therefore,

$$\psi = a\gamma^{\,1}_{\,11} + c_1\varphi + c_2 \ . \tag{a12}$$

Substituting (a6) and (a12) into the last equation of (a3), we obtain

$$\varphi = \left(\gamma^{\,1}_{\,11}\right)^2 + \frac{c_2}{a}\cdot\gamma^{\,1}_{\,11} \ . \tag{a13}$$

Furthermore, by considering the fourth equation of (a3) in light of (a13), we get (a10) (in view of $\left(\gamma^{\,1}_{\,11}\right)' \neq 0$). Therefore, in light of (a10), (a12), and (a13), we get

$$\left.\begin{aligned} \varphi &= \frac{1}{4}\left(\eta^2(z) - \left(\frac{c_2}{a}\right)^2\right) , \\ \psi &= \frac{c_1}{4}\left(\eta^2(z) - \left(\frac{c_2}{a}\right)^2\right) + \frac{a}{2}\left(\eta(z) + \frac{c_2}{a}\right) , \end{aligned}\right\} , \tag{a14}$$

and we get (a11) (cf. (a5), (a6), (a10), and (a14)). □

Now we can prove the following lemma.

**LEMMA A.** *An evolution equation of type (a1) admits a standard connection (with the coefficients* $g^1_{01}(z,p_1)$, $\gamma^1_{00}(z,p_1)$, $\gamma^1_{01}(z)$, *and* $\gamma^1_{11}(z)$, *where* $\frac{d\gamma^1_{11}}{dz} \neq 0$*) if and only if the equation (a1) is of the form (Eq1). Moreover, the connection coefficients are defined by formulas (a10), (a11).*

***PROOF:***

**1.** Let the equation (a1) admit a standard connection for which $\frac{d\gamma^1_{11}}{dz} \neq 0$. Then, as follows from Remark A3, $\gamma^1_{11}(z)$ is defined by (a10), and $\eta'(z) \neq 0$. Substituting (a10) into the second equation of (a3), we get

$$\xi(z) = \frac{\eta''(z)}{\eta'(z)} .$$

So, we see that our evolution equation (a1) is of the form (Eq1). The coefficients $g^1_{01}$, $\gamma^1_{00}$, $\gamma^1_{01}$ are defined by (a11).

**2.** Conversely, suppose we are given an evolution equation of the form (Eq1). Consider the system (a2) together with the condition $\xi = \frac{\eta''}{\eta'}$. It is easy to verify that there is a solution of the system (a2) defined by the formulas (a10), (a11). So (see Theorem 2), an evolution equation (Eq1) admits a standard connection.
This completes the proof of Lemma A. □

**APPENDIX B. Standard connections corresponding to the evolution equation (29) and satisfying** $\frac{d\gamma^1_{11}}{dz} = 0$.

Suppose that a standard connection corresponding to a given evolution equation of type (29) exists, and, moreover, $\frac{d\gamma^1_{11}}{dz} = 0$. Notice that in this case $\frac{d\gamma^1_{01}}{dz} \neq 0$ (see ($25'$)). Here, it will be more convenient to write the equation (29) as follows (we denote $\xi$ by $F'$)

$$z_t - z_{xx} - F'(z)\cdot(z_x)^2 - \eta(z)\cdot z_x - \zeta(z) = 0 . \tag{b1}$$

Notice also that

$$\gamma^1_{11} = a = cons, \tag{b2}$$

as $\frac{d\gamma^1_{11}}{dz} = 0$.

The system (30) (cf. Theorem 2) corresponding to the equation (b1) and satisfying (b2) can be written as follows:

$$\left.\begin{aligned}
&\frac{\partial g_{01}^{1}}{\partial p_1}=0\,,\\
&\frac{\partial g_{01}^{1}}{\partial z}=0\,,\\
&\left(\gamma_{01}^{1}\right)'=\frac{\partial\gamma_{00}^{1}}{\partial p_1}\,,\\
&\left(F'\cdot(p_1)^2+\eta\cdot p_1+\zeta\right)\cdot\left(\gamma_{01}^{1}\right)'=\frac{\partial\gamma_{00}^{1}}{\partial z}\cdot p_1+a\gamma_{00}^{1}-\gamma_{01}^{1}\cdot g_{01}^{1}
\end{aligned}\right\};\qquad\text{(b3)}$$

here $\frac{d\gamma_{01}^{1}}{dz}$ is denoted by $\left(\gamma_{01}^{1}\right)'$.

It follows from the first and the second equations of this system that

$$g_{01}^{1}=b=const;\qquad\text{(b4)}$$

and from the third equation of (b3) that

$$\gamma_{00}^{1}=\left(\gamma_{01}^{1}\right)'\cdot p_1+\psi(z)\,,\qquad\text{(b5)}$$

where $\psi(z)$ is a certain function.

Substituting (b4) and (b5) into the fourth equation of (b3), we get

$$\left(F'\cdot(p_1)^2+\eta\cdot p_1+\zeta\right)\cdot\left(\gamma_{01}^{1}\right)'=\left(\left(\gamma_{01}^{1}\right)''\cdot p_1+\psi'\right)\cdot p_1+a\left(\left(\gamma_{01}^{1}\right)'\cdot p_1+\psi\right)-b\gamma_{01}^{1}\,.\qquad\text{(b6)}$$

Both sides in this relation are polynomials in $p_1$. Compare the coefficients on both sides (remember that $\left(\gamma_{01}^{1}\right)'\neq 0$)

$$\left.\begin{aligned}
&\frac{\left(\gamma_{01}^{1}\right)''}{\left(\gamma_{01}^{1}\right)'}=F'\,,\\
&(\eta-a)\cdot\left(\gamma_{01}^{1}\right)'=\psi'\,,\\
&\zeta\cdot\left(\gamma_{01}^{1}\right)'=a\psi-b\gamma_{01}^{1}
\end{aligned}\right\}.\qquad\text{(b7)}$$

**Remark B1.** *It follows from (b7) that*

$$\gamma_{01}^{1}=c\left(\int e^{F(z)}dz+c_1\right),\qquad\text{(b8)}$$

$$\gamma_{00}^{1}=c\left[e^{F(z)}\cdot p_1+\left(\int(\eta(z)-a)e^{F(z)}dz+c_2\right)\right],\qquad\text{(b9)}$$

*where $c, c_1, c_2$ are constants ($c\neq 0$).*

***PROOF:*** It follows from the first equation of (b7) that

$$\left(\gamma_{01}^{1}\right)'=ce^{F(z)}\,,$$

where $c\neq 0$ (remember that $\left(\gamma_{01}^{1}\right)'\neq 0$). Therefore, (b8) holds.

Substituting (b8) into the second equation of (b7), we get

$$\psi'=c(\eta(z)-a)e^{F(z)}\,,$$

and, hence,

$$\psi(z) = c\left(\int (\eta(z) - a) e^{F(z)} dz + c_2\right). \qquad \text{(b10)}$$

Using (b5), (b8), and (b10), we obtain (b9). □

Now we are ready to prove the following lemma.

**LEMMA B.** *An evolution equation of type (29) admits a standard connection (with the coefficients* $g^1_{01}(z, p_1)$, $\gamma^1_{00}(z, p_1)$, $\gamma^1_{01}(z)$, $\gamma^1_{11}(z)$, *where* $\frac{d\gamma^1_{11}}{dz} = 0$*) if and only if the equation (29) is of the form (Eq2). Moreover, the connection coefficients are defined by the formulas*

$$\left.\begin{aligned} &\gamma^1_{11} = A, \\ &\gamma^1_{01} = c\left(\int e^{F(z)} dz + c_1\right), \\ &\gamma^1_{00} = ce^{F(z)} \cdot p_1 + c\left(\int (\eta(z) - A) e^{F(z)} dz + c_2\right), \\ &g^1_{01} = -A^2 - B \end{aligned}\right\}; \qquad \text{(b11)}$$

here $c, c_1, c_2$ are constants ($c \neq 0$), and $c_1, c_2$ satisfy $(A^2 + B)c_1 + Ac_2 = C$.

***PROOF:***

**1.** Let an equation of type (29) admits a standard connection for which $\frac{d\gamma^1_{11}}{dz} = 0$.

Then, by Remark B1, $\gamma^1_{01}(z)$ and $\gamma^1_{11}(z)$ are of the form (b8) and (b9) respectively. Substituting (b8) and (b10) into the last equation of (b7), we get (in view of $\tilde{n} \neq 0$)

$$\zeta(z) = e^{-F(z)}\left(\int (A\eta(z) + B) e^{F(z)} dz + C\right),$$

where

$$\left.\begin{aligned} &A = a, \\ &B = -a^2 - b, \\ &C = ac_2 - bc_1 \end{aligned}\right\}. \qquad \text{(b12)}$$

So, we see that our evolution equation is of the form (Eq2).

Moreover (see (b12)),

$$a = A, \qquad b = -A^2 - B, \qquad \text{(b13)}$$

and from the last equation of (b12) we get

$$(A^2 + B)c_1 + Ac_2 = C.$$

Formulas (b11) follow from (b2), (b4), (b8), (b9), and (b13).

**2.** Conversely, suppose we are given an evolution equation of type (Eq2). Consider the PDE system (30) under the condition

$$M = F'(z) \cdot (p_1)^2 + \eta(z) \cdot p_1 + e^{-F(z)}\left(\int (A\eta(z) + B) e^{F(z)} dz + C\right).$$

It is easy to verify that there is a solution of this PDE system that is defined by formulas (b11). So, as follows from Theorem 2, an evolution equation of type (Eq2) admits a standard connection.This completes the proof of Lemma B. □

Moscow State M.V.Lomonosov University
Faculty of Mechanics and Mathematics E-mail: arybnikov@mail.ru